\documentclass[letterpaper,12pt]{amsart}
\usepackage[english]{babel}
\usepackage{amsthm}
\usepackage{amsmath}
\usepackage{amssymb}
\usepackage[latin1]{inputenc}
\def\R{\text{$\mathbb{R}$}}

\def\lra{\longrightarrow}

\def\lra{\longrightarrow}

\def\d1#1#2{\frac{d#1}{d#2}}

\def\p1#1#2{\frac{\partial #1}{\partial #2}}
\def\to{t_o}

\def\part{a=\to \leq t_1 \leq \ldots \leq t_{n-1} \leq t_{n}=b}

\def\C{\text{$\mathbb{C}$}}

\newcommand\mf[1]{\mathfrak{#1}}

\newcommand\G{{\rm G}}

\def\o{\omega}

\title{On the moment map on symplectic manifolds}
\author[L. Biliotti]{Leonardo Biliotti}
\address{Dipartimento di Matematica, Universit\`a Politecnica delle
Marche, Via Brecce Bianche, 60131, Ancona  Italy}
\email{biliotti@dipmat.univpm.it}
\begin{document}
\newtheorem{thm}{Theorem}[section]
\newtheorem{prop}[thm]{Proposition}
\newtheorem{lemma}[thm]{Lemma}
\newtheorem{cor}[thm]{Corollary}
\theoremstyle{definition}
\newtheorem{defini}[thm]{Definition}
\newtheorem{notation}[thm]{Notation}
\newtheorem{exe}[thm]{Example}
\newtheorem{conj}[thm]{Conjecture}
\newtheorem{prob}[thm]{Problem}
\theoremstyle{remark}
\newtheorem{rem}[thm]{Remark}
\begin{abstract}
We consider a connected symplectic manifold $M$ acted on by a
connected Lie group $G$ in a Hamiltonian fashion. If $G$ is
compact, we prove give an Equivalence Theorem for the symplectic
manifolds whose squared moment map $\parallel \mu
\parallel^2$ is constant. This result works also in the almost-K\"ahler
setting. Then we study the case when $G$ is a non compact Lie
group acting properly on $M$ and we prove a splitting results for
symplectic manifolds.
\end{abstract}
\maketitle
\section{introduction}
We shall consider symplectic manifolds $(M,\omega)$ acted on by a
connected Lie group $G$ of symplectomorphism.
Throughout this paper we shall assume that the $G$-action is proper and
Hamiltonian, i.e. there exists a moment map $\mu:M \lra \mf g^*$, where
$\mf g$ is the Lie algebra of $G$. In general the matter of
existence/uniqueness of $\mu$ is delicate. However, if $\mf g$
is semisimple, there is a unique moment map  (see \cite{gs}).
If $(M,\omega)$ is a compact K\"ahler manifold and $G$
is a connected compact Lie group of holomorphic isometries,
then the existence problem is resolved
(see \cite{HW}): a moment map exists if and only if $G$ acts
trivially on the Albanese torus Alb($M$) or equivalently every
vector field from $\mf z$, where $\mf z$ is
the Lie algebra of the center of $G$, vanishes at some point in $M$.

If $G$ is compact, we fix an Ad($G$)-invariant
scalar product $\langle \cdot, \cdot \rangle$ on $\mf g$ and we
identify $\mf g^*$ with $\mf g$ by means of $\langle \cdot, \cdot
\rangle$. Then we can think of $\mu$ as a $\mf g -$valued map and it is
natural study the smooth function $f=\parallel \mu
\parallel^2$ which has been extensively used in \cite{ki} to
obtain strong information on the topology of the manifold.

Firstly, we  investigate the symplectic manifolds whose squared moment map
is constant proving the the following Equivalence Theorem.
\begin{thm} \textbf{(Equivalence Theorem)}. \label{eccolo}
Suppose $(M,\omega)$ is a connected symplectic $G$-Ha\-mil\-to\-nian
manifold, where $G$ is a connected compact Lie group acting
effectively on $M$, with moment map $\mu:M \lra \mf g$. Then the
following conditions are equivalent.
\begin{enumerate}
\item $G$ is semisimple and $M$ is $G$-equivariantly symplectomorphic to a
product of a flag manifold and a symplectic manifold which is
acted on trivially by $G$;
\item the squared moment map
$f=\parallel \mu \parallel^2$ is constant;
\item  $M$ is mapped by the moment map $\mu$ to a single coadjoint
orbit;
\item all principal $G$-orbits are symplectic;
\item all $G$-orbits are symplectic.
\end{enumerate}
Moreover, given a $G$-invariant $\omega$-compatible almost
complex structure on $M$, the symplectomorphism in $(1)$
turns out to be an isometry
with respect to the induced Riemannian metric
while $(4)$ and $(5)$ become:
all $G$-orbits (resp. principal $G$-orbits) are complex.
\end{thm}
In order to prove the above theorem, we need the following result, which
might have an indipendent interest.
\begin{prop} \label{max}
Let $(M,\omega)$ be a symplectic manifold and let $G$ be connected
compact Lie group acting in a Hamiltonian fashion on $M$ with
moment map $\mu$. Assume $x \in M$ realizes a local maximum of the
smooth function $f=\parallel \mu \parallel^2$. Then $G \cdot x$ is
symplectic and there exists a neighborhood $Y_o$ of $x$ such that
$G\cdot ( Y_o \cap \mu^{-1}(\mu(x)))$ is a symplectic submanifold
which is $G$-equivariantly symplectomorphic to a flag manifold and
a symplectic manifold which is acted on trivially by $G$.
Moreover, if $x\in M$ realizes the maximum of $f=\parallel \mu
\parallel^2$ or any $z\in \mu^{-1}(\mu(x))$ realizes a local
maximum of $f=\parallel \mu
\parallel^2$, then
\begin{enumerate}
\item $\mu^{-1}(\mu(x))$ is a  symplectic submanifold of $M$;
\item $G\cdot \mu^{-1}(\mu(x))$ is a symplectic submanifold of $M$ which is
$G-$equiva\-ri\-an\-tly symplectomorphic to $(Gx
\times\mu^{-1}(\mu(x)), \omega_{|_{G \cdot x}}+
\omega_{|_{\mu^{-1}(\mu(x))}})$.
\end{enumerate}
\end{prop}
These results generalize ones given in \cite{gp} and \cite{gb}.

One may try to prove Proposition \ref{max} assuming only Ad$(G)$
is compact; this means $G$ is covered for a compact Lie group and
a vector group which lies in the center (see \cite{CE}). However,
if $G$ acts properly on $M$, then the existence of a
$G$-symplectic orbit implies that $G$ must be compact. Indeed if
$G\cdot x=G/G_x$ is symplectic, from Proposition \ref{KKS}, then $( G_x
)^o$, the connected component of stabilizer group of $x$ which contains the
identity, coincides with $(G_{\mu(x)})^o$. Since $(Z(G))^o \subset
(G_{\mu(x)})^o$, we conclude that $G$ is compact.

Then we study the case when $G$ is a
non compact Lie group acting properly and in a Hamiltonian fashion
on $M$ and we prove the following result. 
\begin{thm} \label{ultimo}
Let $(M,\omega)$ be a symplectic manifold and let $G$ be a connected
non compact Lie group acting properly and in a Hamiltonian fashion on $M$ with
moment map $\mu$. Assume also $G\cdot \alpha$ is a locally closed coadjoint
orbit for every $\alpha \in \mf g^*$. Then the following
conditions are equivalent.
\begin{enumerate}
\item All $G$-orbits are symplectic;
\item all principal $G$-orbits are symplectic;
\item $M$ is mapped by the moment map $\mu$ to a single coadjoint orbit;
\item let $x$ be a regular point of $M$.
Then $G \cdot x$ is a symplectic orbit, $\mu^{-1}(\mu(x))$
is a symplectic submanifold on which $G_x$ acts trivially
and the following $G$-equivariant application
\[
\phi: G\cdot x \times \mu^{-1}(\mu(x)) \lra M,\
\phi([gx,z])=gz,
\]
is surjective and satisfies
\[
\phi^* (\omega)= \omega_{|_{G \cdot x}} + \omega_{|_{\mu^{-1}(\mu(x))}}.
\]
\end{enumerate}
If $G$ is a reductive Lie group acting effectively on
$M$, then $G$ has to be semisimple and $\phi$ becomes a $G$-equivariant
symplectomorphism.
Moreover, if $N(G_x)/G_x$ is finite,
given a $G$-invariant $\omega$-compatible almost
complex structure on $M$, the symplectomorphism in $(4)$
turns out to be an isometry
with respect to the induced Riemannian metric
while $(1)$ and $(2)$ become:
all $G$-orbits (resp. principal $G$-orbits) are complex.
\end{thm}
The assumption $G\cdot \alpha$ is a locally
closed coadjoint orbit is needed to applying the symplectic slice
and the symplectic stratification of the reduced spaces given in
\cite{lb}. Observe that the condition of a coadjoint orbit being
locally closed is automatic for reductive group and for their
product with vector spaces. There exists an example of a solvable
group due to Mautner \cite{mau} p.$512$, with non-locally closed
coadjoint orbits.

Finally, as an immediate corollary of
Theorem \ref{eccolo} and
Theorem \ref{ultimo}, we give the following  splitting result.

Let $G$ be a non
compact semisimple Lie group. The Killing form $B$ on $\mf g$ is a
non-degenerate Ad$(G)$-invariant bilinear form. Therefore, we may identify
$\mf g$ with $\mf g^*$ by means of $-B$ and we may think $\mu$ as
a $\mf g$-valued map.
\begin{cor} \label{cccp} Let $M$ be a symplectic manifold acted
on by connected non compact semisimple Lie group $G$, properly and in a
Hamiltonian fashion with moment map $\mu$. If $f=\parallel \mu
\parallel^2$ is constant and any element which lies in the image  of
the moment map $\mu$ is elliptic, then all $G$-orbits are symplectic and $M$
is $G$-equivariantly symplectomorphic to a product of a flag manifold and a
symplectic manifold which is acted on trivially by $G$.
Moreover, if $N(G_x)/G_x$ is finite,
given a $G$-invariant $\omega$-compatible almost
complex structure on $M$, the application $\phi$ turn out to be an isometry
with respect to the induced Riemannian metric and
all $G$-orbits are complex.
\end{cor}
\section{Proof of the main results} Let $M$ be a connected
differential manifold equipped with a non-degenerate closed
$2-$form $\o$. The pair $(M,\o)$ is called \emph{symplectic manifold}.
Here we consider a finite-dimensional connected Lie group acting
smoothly and properly on $M$ so that $g^* \omega=\omega$ for all
$g \in G,$ i.e. $G$ acts as a group of canonical or symplectic
diffeomorphism.

The $G$-action is called \emph{Hamiltonian}, and we said that $G$
acts in a Hamiltonian fashion on $M$ or $M$ is $G$-Hamiltonian,
if there exists a map
$$
\mu: M \lra \mf{g}^*,
$$
which is called moment map, satisfying:
\begin{enumerate}
\item for each $X\in \mf{g}$ let
\begin{itemize}
\item $\mu^{X}: M \lra \R,\ \mu^{X} (p)= \mu(p) (X),$ the
component of $\mu$ along $X,$ \item $X^\#$ be the vector field on
$M$ generated by the one para\-me\-ter subgroup $\{ \exp (tX):
t \in \R \} \subseteq G$.
\end{itemize}
Then
$$
{\rm d} \mu^{X}= {\rm i}_{X^\#} \o ,
$$
i.e. $\mu^{X}$ is a Hamiltonian function for the vector field
$X^\#.$
\item $\mu$ is $G-$equivariant, i.e. $\mu (gp)=Ad^* (g)
(\mu(p)),$ where $Ad^*$ is the coadjoint representation on $\mf g^*$.
\end{enumerate}

Let $x \in M$ and ${\rm d} \mu_x: T_x M \lra T_{\mu(x)} \mf g^*$
be the differential of $\mu$ at $x$. Then
\[
Ker {\rm d}\mu_x= (T_x G \cdot x )^{\perp_{\omega}}:= \{ v \in T_x M:
\omega (v,w)=0,\ \forall w \in  T_x G \cdot x\}.
\]
If we restrict $\mu$ to a $G-$orbit $G \cdot x$, then we have the following
homogeneous fibration
\[
\mu: G \cdot x \lra Ad^*(G)\cdot \mu(x)
\]
and the restriction of the ambient symplectic form $\omega$ to the orbit
$G \cdot x$ equals the pullback by the moment map $\mu$ of the symplectic
form on the coadjoint orbit through $\mu(x)$:
\begin{equation} \label{KKS}
\omega_{|_{G \cdot x}}= \mu^* (\omega_{ Ad^*(G) \cdot \mu(x)})_{|_{G \cdot x}},
\end{equation}
see \cite{lb} p. 211, where $\omega_{G \cdot \mu(x)}$ is the
Kirillov-Konstant-Souriau (KKS) symplectic form on the coadjoint
orbit of $\mu(x)$ in $\mf g^*$. This implies the following well-known fact,
see \cite{gs}.
\begin{prop} \label{c1}
The orbit of $G$ through  $x\in M$ is symplectic if and only if
the stabilizer group of $x$ is an open subgroup of the stabilizer
of $\mu(x)$ if and only if
the moment map restricted to $G \cdot x$ into $G \cdot \mu(x)$,
$\mu_{|G \cdot x}:G \cdot x \longrightarrow G \cdot \mu(x)$,
is a covering map.
In particular if $G$ is compact or semisimple, then $G_x=G_{\mu(x)}$ so this
implies that $\mu_{|_{G \cdot x}}:G \cdot x \longrightarrow G \cdot \mu(x)$
is a diffeomorphism.
\end{prop}
\begin{proof}
The first affirmation follows immediately from (\ref{KKS}). If $G$
is compact or semisimple, then $G_{\mu(x)}$ is connected so this implies
that the two stabilizer groups are the same.
\end{proof}
We now give the proof of Proposition \ref{max} \\
$\ $ \\
\emph{Proof of Proposition \ref{max}}. Let $\beta=\mu(x)$ and let
$G_x$ be the isotropy group at $x$. From the local normal for the
moment map, see  \cite{lb}, \cite{gs}, \cite{or} and \cite{sl},
there exists a neighborhood of the orbit $G \cdot x$ which is
e\-qui\-va\-riantly symplectomorphic to a  neighborhood $Y_o$ of
the zero section of $ (Y=G \times_{G_x}( \mathbf{\mf q} \oplus V),
\tau)$ with the $G$-moment map $\mu$ given by
\[
\mu([g,m,v])=Ad(g)(\beta+m+ \mu_{V}(v)),
\]
where $\mathbf{\mf q }$ is a summand in the $G_x$-equivariant splitting
$\mf g=\mf g_{\beta} \oplus \mathbf{\mf s}$ $= \mf g_x \oplus \mathbf{\mf q}
\oplus \mathbf{\mf s}$ and  $\mu_V$ is the moment map of the
$G_x$-action on the symplectic subspace $V$ of
$((T_x G\cdot x)^{\perp_{\omega}},\omega(x))$.
Note that $V$ is isomorphic to the quotient $((T_x
G \cdot x)^{\perp_{\omega}} /
((T_x G \cdot x)^{\perp_{\omega}} \cap T_x G \cdot x))$.

In the sequel we denote by $\omega_V= \omega(x)_{|_{V}}$ and
shrinking $Y_o$ if
necessary, we may suppose that $[e,0,0]$ is a maximum of the
smooth function $f=\parallel \mu \parallel^2$ in $Y_o$.

We first want to prove that $G \cdot x$ is symplectic. Then we shall prove
$\mathbf{\mf q}= \{ 0 \}$.

Let $m \in \mathbf{\mf q}-\{0\}$.
Then for every $\lambda \in \R$ we have
\[
f(e,\lambda m,0)=
\parallel \beta \parallel^2 + \lambda^2 \parallel m \parallel^2 +
\lambda \langle m,\beta \rangle \leq \parallel \beta \parallel^2
\]
and therefore
\[
\lambda^2 \parallel m \parallel^2 + \lambda \langle m, \beta
\rangle \leq 0,
\]
for every $\lambda \in \R$ which is a contradiction.
Hence $G \cdot x$ is symplectic and by Proposition
\ref{c1} $G_x=G_{\beta}.$

Note that any $y\in Y^{\beta}_o= Y_o \cap \mu^{-1}(\beta)$
is a local maximum for $f$.
Then $G_y=G_x$ for every $y\in Y^{\beta}_o$, i.e. $G\cdot y$ is symplectic, and
a $G-$orbit through an element of $Y^{\beta}_o$ intersects
$\mu^{-1}(\beta)$ in at most one point. Indeed, if both $x\in Y^{\beta}_o$ and
$kx$ lie in $\mu^{-1}(\beta)$, then, by the $G$-equivariance of $\mu$,
we have  $\mu(kx)=\beta=k\mu(x)=k\beta$, proving $k \in G_x$.
From this it follows that the map
\[
\phi:G\cdot x \times Y^{\beta}_o \longrightarrow G\cdot Y^{\beta}_o
\]
is well-defined and bijective.

Let us now return to the local normal form.
Let $y\in Y^{\beta}_o$. We know that $G_y=G_x$ and $G\cdot y$ is symplectic.
Hence there exists a neighborhood of $G\cdot y$ which is $G$-invariant
symplectomorphic to a neighborhood $Y'$ of the zero section of
$(Y=G \times_{G_x} V, \tau)$ with the $G$-moment map
given by
\[
\mu([g,v])=Ad(g)(\beta + \mu_V (v)).
\]
Shrinking $Y'$, if necessary, we may assume $Y' \subseteq Y_o$ and,
see Proposition 13 in \cite{lb} p.216, the intersection
of the set $\mu^{-1}(G\cdot \beta)$ with $Y'$ is of the form
\[
G \cdot \mu^{-1}(\beta) \cap Y'=\{[g,v]\in Y_o:\mu_V (v)=0 \}.
\]
Let $Y^{(G_x)}=\{m\in Y:\ (G_m )=(G_x ) \}$, i.e. $G_m$ is $G$-conjugate to
$G_x$.
It is easy to check that
\begin{equation} \label{ote}
Y^{(G_x)}=G \times_{G_x} V^{G_x} \cong G/G_x \times V^{G_x},
\end{equation}
where $V^{G_x}=\{x\in V:G_m=G_x \}$, and  $\mu(Y^{( G_x )})=G \cdot \beta$.
Therefore,
since $Y^{\beta}_o \subseteq M^{G_x}$, we have
\begin{equation} \label{splitting}
Y' \cap \mu^{-1}(G \cdot \beta)= Y^{(G_x)} \cap Y'=(Y' )^{(G_x )}
\end{equation}
and
\[
Y' \cap \mu^{-1}(\beta)=Y' \cap V^{G_x}.
\]
This implies both $Y^{\beta}_o$ and $G \cdot Y^{\beta}_o$
are symplectic submanifolds of $M$. Indeed, from
the above discussion we conclude that
$T_y Y^{\beta}_o=V^{G_x}$ and
the tangent space at $y$ of
$G\cdot Y^{\beta}_o$ splits as
\[
T_y G \cdot Y^{\beta}_o = T_y G\cdot y
\stackrel{\perp_{\omega}}{\oplus}
T_y Y^{\beta}_o.
\]
Here we have used that
$T_y Y^{\beta}_o \subset ( T_y G \cdot y )^{\perp_{\omega}}= Ker \mathrm{d}\mu_y$
and  $G \cdot y$ is symplectic.

Since
\begin{equation} \label{splitting2}
\tau_{|_{G/G_x \times_{G_x} V^{G_x }}}= \omega_{|_{G \cdot x}} +
(\omega_{V})_{|_{ V^{ G_x } }},
\end{equation}
see Corollary $14$ p. 217 \cite{lb}, from (\ref{KKS}), (\ref{ote}),
(\ref{splitting}),
and (\ref{splitting2})
we obtain that $\phi$ is a symplectomorphism.

Now assume that $x\in M$ realizes the maximum of $f$ or any
$z\in \mu^{-1}(\mu(x))$ is a local maximum of $f$. Let
$\beta=\mu(x)$. Using the same argument as before,
we may prove
$G \cdot z$ is symplectic, $G_z = G_x =G_{\beta}$ for every
$z \in \mu^{-1}(\beta)$ and a $G$-orbit intersects
$\mu^{-1}(\beta)$ in at most one point.
It follows that the following application
\[
\phi:G \cdot x \times \mu^{-1}(\beta) \lra G \cdot \mu^{-1}(\beta),\ \
\phi(gG_x,z)=gz
\]
is a $G$-equivariant diffeomorphism. We shall prove that $\phi$ is a
symplectomorphism.

The set $ \mu^{-1}( G \beta )\cap M^{(G_x ) }$ is a manifold
of constant rank and the quotient
\[
(M_{\beta} )^{(G_x)}:=(
\mu^{-1}( G \cdot \beta )\cap M^{(G_x ) } )/G,
\]
is a symplectic manifold, see Corollary $14$ in \cite{lb}.
Since $\mu^{-1}(\beta) \subset M^{G_x }$, we have
\[
\G \cdot \mu^{-1}(\beta)=\mu^{-1}(G \cdot \beta )\cap M^{(G_x ) },
\]
i.e.  $G\cdot \mu^{-1}(\beta)$ is a submanifold,
and finally $\beta$ is a regular value
of
\[
\mu_{|_{ \mu^{-1} ( G \cdot \beta )}} : G\cdot \mu^{-1}(\beta)
\longrightarrow G\cdot \beta.
\]
Therefore
$\mu^{-1}(\beta)$ is a submanifold of $M$ and for
every $z\in \mu^{-1}(\beta)$ the tangent space of $G\cdot \mu^{-1}(\beta )$
splits as
\begin{equation} \label{split}
T_z G \cdot z  \stackrel{\perp_{\omega}}{\oplus}
T_z \mu^{-1}(\beta) =T_z G \cdot \mu^{-1}(\beta).
\end{equation}
Since $T_z \mu^{-1}(\beta)=(V)^{G_x}$ and $G \cdot z$ is symplectic, we
conclude that both $G \cdot \mu^{-1}(\beta)$ and $\mu^{-1}(\beta)$ are
symplectic submanifolds of $M$. Moreover, from (\ref{KKS})  and
(\ref{split}) we obtain that $\phi$ is a $G$-equivariant symplectomorphism.
\begin{flushright}
$\square$
\end{flushright}
\emph{Proof of Theorem \ref{eccolo}}. ((1)$\iff$(2)).
(1)$\Rightarrow$(2)
is trivial. \\
(2)$\Rightarrow$(1).
Assume that the square of the moment map is constant.
Let $x\in M$. By the argument used in the proof of
Proposition \ref{max}, we have $G \cdot x$ is symplectic  and
$G_x=G_{\mu(x)}$. Therefore all $G$-orbits are symplectic
((2)$\Rightarrow$(5)) and
the center acts trivially on $M$, i.e. $G$ is semisimple. Indeed,
coadjoint orbits are of the form $G/C(T)$, where
$C(T)$ is the centralizer of the torus $T$. In particular
$Z(G)\subset G_x$ for every $x\in M$.

We want to show that the manifold $M$ is mapped by the moment map to a single
coadjoint orbit ((2)$\Rightarrow$(3)).

Let $G \cdot x$ be a principal orbit. Since $G_x$ acts trivially
on the slice, from the local normal form for the moment map, in a $G$-invariant
neighborhood of $G\cdot x$ the moment map is given by
\[
\mu([g,v])=Ad(g)(\beta).
\]
This proves that  there exists a $G$-invariant neighborhood of $G \cdot x$
which is mapped to a single coadjoint orbit. It is well-known that the set
$M^{(G_x)}$ is an open dense subset of $M$ and
$M^{(G_x)}/G$ is connected (\cite{Path}). Since $\mu$ is
$G$-equivariant, it induces a continuous application
\[
\overline\mu: M^{(G_x)}/G \lra \mf g/G,
\]
which is locally constant. Hence $\overline\mu(M^{( G_x )} /G)$
is constant so $\overline\mu(M /G )$ is.
Thus $M$ is mapped by $\mu$ to a single
coadjoint orbit; in particular $M=G\cdot \mu^{-1}(\beta)$.
Note that this argument proves (4)$\Rightarrow$(3).

Let $x\in M$.  As in the proof of Proposition \ref{max}, from (\ref{KKS}),
(\ref{splitting}),
(\ref{splitting2}) and (\ref{split}), the following
application
\[
\phi: G \cdot x \times \mu^{-1}(\mu(x)) \lra  M, \ \ \ (gx , z)
\lra gz,
\]
is the desired $G$-equivariant symplectomorphism. \\
(2)$\Rightarrow$(3), (2)$\Rightarrow$(5) and (4)$\Rightarrow$(3)
follow from the above discussion
while (3)$\Rightarrow$(2) and (5)$\Rightarrow$(4) are easy to check.\\
(3)$\Rightarrow$(5). In the sequel we follow
the notation introduced in the proof of
the Proposition \ref{max}. Let $G\cdot x$ be a $G$-orbit and let $Y'$ be the
neighborhood of the zero section in
$(Y=G\times_{G_x} (\mathbf{\mf q}\oplus V,\tau)$ which
is $G$-equivariant symplectomorphic to a neighborhood of $G\cdot x$.
The moment map $\mu$ in $Y'$ is given by the formula
 \[
\mu([g,m,v])=Ad(g)(\beta+m+ \mu_{V}(v)).
\]
From Proposition 13 in \cite{lb}, shrinking $Y'$ if necessary, we have
\[
\mu^{-1}(G\cdot \beta) \cap Y'=\{[g,m,v]:m=0 \ \mathrm{and}\ \mu_V(v)=0\}.
\]
Since $M$ is mapped by the moment map $\mu$ to a single coadjoint orbit
$G\cdot \beta$, we conclude that $\mathbf{\mf q}=\{0\}$ and therefore
$G\cdot x$ is symplectic.

Now assume that $M$ is a K\"ahler manifold.
Then $\omega=g( J\cdot , \cdot)$ where $J$ is
the integrable complex structure.

We shall prove that if $x\in M$, then
the application
\[
\phi:Gx \times \mu^{-1}(\mu(x)),\ \ \phi(gx,z)=gz,
\]
is an isometry.

From the argument used in the proof of
Proposition \ref{max}, we have
that for every
$y\in \mu^{-1}(\mu(x))$, $T_{y} \mu^{-1}(\mu(x))=T_y M^{G_x }$. Indeed,
$G_y=G_x=G_{\mu(x)}$ and $G_{\mu(x)}$ centralizes a torus; therefore
$N(G_y)/G_y$ is finite. Since the complex totally geodesic submanifold
$M^{G_x }$ is given by
\[
M^{G_x}=N(G_x )/G_x \times V^{G_x}=N(G_x)/G_x \times \mu^{-1}(\mu(x)),
\]
where the last equality follows from the fact that $\phi$ is a
$G$-equivariant symplecomorphism,
we conclude that  the connected component of
$\mu^{-1}(\mu(x))$ which contains $x$ is the connected component of
$M^{G_x}$ which contains $x$. Therefore $\mu^{-1}(\mu(x))$ is a complex
totally geodesic submanifold of $M$.

Now, we show that all $G$-orbits are complex.

Let $V\in T_{y} \mu^{-1}(\mu(x))$ and let $X^\#$ be a tangent vector to
$T_y G \cdot y$ induced  by $X\in \mf g$. Then
\begin{equation} \label{splic}
0={\rm d}\mu_y (V)(X)=\omega(X^\#,V)=g(J(X^\# ),V)=g(X^\#,J(V))
\end{equation}
which implies that $T_y G \cdot y=T_y \mu^{-1}(\mu(x))^{\perp}$ so
$G \cdot y$ is complex.

Finally we prove that $\phi$ is an isometry. Since $\phi$ is $G$-equivariant
and $G$ acts by isometries, it is enough to prove that d$\phi_{(x,z)}$ is an
isometry for every $z\in \mu^{-1}(\mu(x))$. Note that the tangent space of $M$
splits as
\begin{equation} \label{splitc}
T_z M = T_z G\cdot z \stackrel{\perp}{\oplus} T_z \mu^{-1}(\mu(x)),
\end{equation}
for every $z\in \mu^{-1}(\mu(x))$.
Hence it is sufficient to prove that the Killing vector
field from $\xi \in \mf g$ has constant norm along $\mu^{-1}(\mu(x))$.

Let $\xi \in \mf g$ and let $X$ be a vector field  tangent to
$\mu^{-1}(\mu(x))$. Then $[\xi^\#,X]=0$ since $\phi$ is a
$G$-equivariant diffeomorphism. Now, given $\eta^\#$ be such that
$J(\eta^\#)=\xi^\#$, by the closeness of $\omega$ we have
\[
0=d\omega(X,\eta^\#, \xi^\#)=Xg(\xi^\#, \xi^\#)
\]
which implies that the Killing field from $\xi \in \mf g$
has constant norm along $\mu^{-1}(\mu(x))$.

In \cite{lb} it was proved that there exists a $G$-invariant almost
complex structure $J$ adapted to $\omega$, i.e.
$\omega(J\cdot, J\cdot) =\omega(\cdot, \cdot)$ and $\omega(\cdot,J\cdot)=g$
is a Riemannian metric. Since $T_y \mu^{-1}(\mu(x))=T_y M^{G_x }$,
$\mu^{-1}(\mu(x))$ is $J$-invariant. This allow us to conclude that the
splitting (\ref{splitc}) holds, hence every $G$-orbit is $J$-invariant.
Now, the fact $\phi$ is an isometry follows as before, proving the
result in the almost-K\"ahler setting.
\begin{flushright}
$\square$
\end{flushright}
\emph{Proof of Theorem \ref{ultimo}}.
(1)$\iff$(2)$\iff$(3)
follow using the same arguments in the proof of
the Theorem \ref{eccolo} while
(4)$\Rightarrow$(3) is easy to check.
We shall prove (3)$\Rightarrow$(4).

Let $x\in M$ be a regular point and
let $\beta=\mu(x)$. As in the proof of the Theorem
\ref{eccolo}, $M=G \cdot \mu^{-1}(\beta)$ and any orbit is symplectic. This
implies that
$\beta$ is a regular value of the application
\[
\mu: M \longrightarrow G\cdot \beta.
\]
Therefore $\mu^{-1}(\beta)$ is a closed submanifold whose tangent
space is given by
\[
T_y \mu^{-1}(\beta)=Ker \mathrm{d}\mu_y =( T_y G \cdot y )^{\perp_{\omega}}
\]
and the tangent space of $M$ splits as
\begin{equation} \label{c2}
 T_y M = T_y G \cdot y
\stackrel{\perp_{\omega}}{\oplus} T_y \mu^{-1}(\beta).
\end{equation}
for every $y\in \mu^{-1}(\beta)$. In particular $\mu^{-1}(\beta)$
is symplectic.

Now, we show that $G_x$ acts trivially on $\mu^{-1}(\beta)$.

Note that $(G_y )^o=(G_x)^o=(G_{\beta})^o$, for every $y\in \mu^{-1}(\beta)$,
due the fact that $G\cdot y$  is symplectic, $G_y \subset G_{\beta}$
since $\mu$ is $G$-equivariant, and
$\mu^{-1}(\beta)$ is connected since both $G$ and
$M=G \cdot \mu^{-1}(\beta)$ are.

From the slice theorem, see \cite{Path},
there exists a neighborhood $U$ of the regular point $x$
such that $(G_z)=(G_x)$ $\forall z\in U$. Assume that we may find
a sequence $x_n \rightarrow x$ in $\mu^{-1}(\beta)$
and a sequence $g_n \in G_{x_n}-G_x \subseteq G_{\beta}$
such that $g_n x_n=x_n$.
Since the $G$-action is proper we may assume that
$g_n \rightarrow g_o$ which lies in $G_x$.
In particular the sequence $g_n$ converges to $g_o$ in
$G_{\beta}$ as well, since it is a closed Lie group.

Now, $G_x$ is an open subset of $G_{\beta}$, since
$(G_x )^o= (G_{\beta})^o$; therefore there exists $n_o$ such that
$g_n \in G_x$ for $n\geq n_o$ which is an absurd.

Thus, there exists an open
subset $U'$ of $x$ in $\mu^{-1}(\beta)$ such that $G_z=G_x$, $\forall z\in U'$.
Since $G_x$ is compact we conclude that $G_x$ acts trivially on
$\mu^{-1}(\beta)$.

It follows that the application
\[
\phi:G \cdot x \times \mu^{-1}(\beta) \longrightarrow M,\ \ \ \phi(gG_x,z)=gz
\]
is well-defined, smooth and $G$-equivariant.
Moreover, from (\ref{KKS}) and (\ref{c2}) we get
that
\begin{equation} \label{pio}
\phi^* (\omega)= \omega_{|_{G \cdot x}} + \omega_{|_{\mu^{-1}(\beta )}}.
\end{equation}
Note also
$Z(G)^o \subseteq (G_{\mu(x)})^o= (G_x )^o$, so $G$ must be semisimple whenever
$G$ is a reductive Lie group acting effectively on $M$ and $\phi$ becomes
a symplectomorphism from Proposition \ref{c1} since a $G$-orbit intersects
$\mu^{-1}(\beta)$ in one point.

Now assume $M$ is almost-K\"ahler and $N(G_x)/G_x$ is finite. The
same arguments used  in the proof of Theorem \ref{eccolo}
show that
$\mu^{-1}(\beta) \cap M^{G_x}$ is complex, $G \cdot y$ is complex for
every $y\in \mu^{-1}(\beta) \cap M^{G_x}$ and the following map
\[
\overline{\phi}=\phi_{|_{G \cdot x \times (\mu^{-1}(\beta) \cap M^{G_x})}}:
G \cdot x \times (\mu^{-1}(\beta) \cap M^{G_x})
\lra M^{(G_x )}, \ \phi([gx,z])=gz,
\]
is $G$-equivariant and it satisfies (\ref{pio}); therefore
a local diffeomorphism. Then we may use the same argument in the proof of
the Theorem \ref{eccolo} to prove that $\overline\phi$ is an isometry.
Since $M^{( G_x )}$ is an open dense subset of $M$,
we obtain that $\mu^{-1}(\beta)$
is complex, so all $G$-orbits are, since the symplectic splitting
(\ref{splitc}) turns out to be $g$-orthogonal, where
$g=\omega(\cdot, J\cdot)$ is the induced Riemannian metric,
and finally we get $\phi$ is an isometry.
\begin{flushright}
$\square$
\end{flushright}
\emph{Proof of Corollary \ref{cccp}} We recall that an
element $X\in \mf g$
is called \emph{elliptic} if $\mathrm{ad}(X)\in \mathrm{End}({\mf g}^{\C})$
is diagonalizable and all eigenvalues are purely immaginary. Then $G\cdot X$
is called an \emph{elliptic orbit}. See \cite{ne} and \cite{ko} for more
details about elliptic orbits.

Let $\mf g= \mf k \oplus \mf p$ be a
Cartan decomposition of the Lie algebra of $G$. Since any elliptic element
is conjugate to an element of $\mf k$, we have that
the squared moment map $f=\parallel \mu \parallel^2$ is positive.
Therefore, using the same argument in the proof of Theorem
\ref{max}, we conclude that all $G$-orbits are
symplectic.
Now the result follows from Theorem \ref{ultimo}.
\begin{flushright}
$\square$
\end{flushright}

\end{document}